\documentclass[11pt,twoside]{article} 

\usepackage{epsfig}
\usepackage[usenames]{color}
\usepackage{amsmath}
\usepackage{amssymb}
\usepackage{mathrsfs}
\usepackage{amsfonts}
\usepackage{amsthm}
\usepackage{lscape}
\usepackage{graphicx}
\usepackage{hyperref}
\hypersetup{colorlinks=false,pdftitle="Reverse entropy power inequalities",pdftex}

\theoremstyle{plain} %default is 'plain' which italicizes body text
\newtheorem{theorem}{Theorem}[section]

\newtheorem{corollary}[theorem]{Corollary}

\theoremstyle{remark} %default is 'plain' which italicizes body text

\newcommand{\bc}{\begin{center}}
\newcommand{\ec}{\end{center}}
\newcommand{\bt}{\begin{tabular}}
\newcommand{\et}{\end{tabular}} 
\newcommand{\bea}{\begin{eqnarray}}
\newcommand{\eea}{\end{eqnarray}}
\newcommand{\bean}{\begin{eqnarray*}}
\newcommand{\eean}{\end{eqnarray*}}

\newcommand{\ba}{\begin{array}}
\newcommand{\ea}{\end{array}}

\def\be{\begin{eqnarray}}
\def\ee{\end{eqnarray}}
\def\ben{\begin{eqnarray*}}
\def\een{\end{eqnarray*}}

\newcommand{\nth}{\frac{1}{n}}

\newcommand{\RL}{{\mathbb R}}

\def\elabel#1{\label{e:#1}}

\def\sq{$\Box$}

\def\qed{\ifmmode\sq\else{\unskip\nobreak\hfil
\penalty50\hskip1em\null\nobreak\hfil\sq
\parfillskip=0pt\finalhyphendemerits=0\endgraf}\fi\par\medbreak}

\newsavebox{\junk}
\savebox{\junk}[1.6mm]{\hbox{$|\!|\!|$}}

\def\til={{\widetilde =}}

\def\half{{\mathchoice{\textstyle \frac{1}{2}}%
{\frac{1}{2}}%
{\hbox{\tiny $\frac{1}{2}$}}%
{\hbox{\tiny $\frac{1}{2}$}} }}

 \def\eq#1/{(\ref{#1})}

\def\eq#1/{(\ref{e:#1})}

\newcommand{\beqn}[1]{\notes{#1}%
\begin{eqnarray} \elabel{#1}}

\newcommand{\eeqn}{\end{eqnarray} }

\newcommand{\beq}[1]{\notes{#1}%
\begin{equation}\elabel{#1}}

\newcommand{\eeq}{\end{equation}} 

\def\bdes{\begin{description}}
\def\edes{\end{description}}

\def\notes#1{}

\newcommand{\calN}{\mathcal{N}}

\renewcommand{\k}{\kappa}
\newcommand{\R}{\RL}

\def\P{{\bf P}}
\def\R{{\bf R}}
\def\k{\kappa}
\newcommand{\en}{\ee}
\def\bee{\begin{eqnarray*}}
\def\ene{\end{eqnarray*}}
\newcommand{\h}{\widetilde{h}}
\def\E{{\bf E}\,}

\numberwithin{equation}{section}

\begin{document}

\title{Dimensional behaviour of entropy and information}
\author{Sergey Bobkov\thanks{S. G. Bobkov is with the
              School of Mathematics, University of Minnesota, 127 Vincent Hall, 
              206 Church St. S.E., Minneapolis, MN 55455 USA.
              Email: {\tt bobkov@math.umn.edu}.}
and Mokshay Madiman\thanks{M. M. Madiman is with the
          Department of Statistics, Yale University, James Dwight Dana House,
          24 Hillhouse Ave, New Haven, CT 06511 USA.
              Tel.: +1-203-432-0639.
              Fax: +1-203-432-0633. 
          Email: {\tt mokshay.madiman@yale.edu}.}
\thanks{S.B. was supported in part by U.S. National Science Foundation grant DMS 0706866, 
and M.M. was supported in part by a Junior Faculty Fellowship from Yale University.} 
}
\date{}
\maketitle

\noindent{\bf Note}: A slightly condensed version of this paper is in press and will appear in 
the  {\it Comptes Rendus de l'Acad\'emies des Sciences Paris, S\'erie I Mathematique}, 2011.

\begin{abstract}
We develop an information-theoretic perspective on some questions in convex geometry,
providing for instance a new equipartition property for
log-concave probability measures, 
some Gaussian comparison results for log-concave measures,
an entropic formulation of the hyperplane conjecture, and
a new reverse entropy power inequality for log-concave measures
analogous to V.~Milman's reverse Brunn-Minkowski inequality.
\end{abstract}

\section{Introduction}
\label{sec:intro}

This Note announces some of the results obtained in \cite{BM10:smb,BM09:maxent,BM10:repi}.
%In order to compactly convey them, we organize them around the following analogue for log-concave
%measures of Milman's reverse Brunn-Minkowski inequality for convex bodies. 
Given a random vector $X$ in $\RL^n$ with density $f(x)$, the entropy power
is defined by
$
\calN(X) = e^{2h(X)/n} ,
$
where, with a common abuse of notation, we write $h(X)$ for  
the Shannon entropy $h(f):=-\int_{\RL^n} f\log f$.

\vspace{.15cm}
\begin{theorem}\label{thm:repi}
If $X$ and $Y$ are independent random vectors 
in $\RL^n$ with log-concave densities, 
there exist affine entropy-preserving %(and volume-preserving) 
maps $u_i:\RL^n \rightarrow \RL^n$ such that 
\ben
\calN\big(\widetilde X + \widetilde Y\big)\, \leq\, C\, (\calN(X) + \calN(Y)),
\een
%then
%\vspace{-.25cm}
%\ben
%\exists \,\,\text{affine entropy-preserving maps } u_i:\RL^n \rightarrow \RL^n \text{ such that }\,
%\calN\big(\widetilde X + \widetilde Y\big)\, \leq\, C\, (\calN(X) + \calN(Y)),
%\een
%\vspace{-.9cm}
where $\widetilde X = u_1(X)$, $\widetilde Y = u_2(Y)$, and where
$C$ is a universal constant.
\end{theorem}
\vspace{.15cm}

Observe that the Shannon--Stam entropy power inequality \cite{Sta59} implies that 
$
\calN\big(\widetilde X + \widetilde Y\big)\,\geq \, \calN(X) + \calN(Y)
$
is always true. Thus Theorem~\ref{thm:repi} may be seen as a reverse
entropy power inequality for log-concave measures. The proof of this
assertion, outlined in Section 3, is based on a series of propositions
introduced in Section 2 including V.~Milman's result on the existence of $M$-ellipsoids. 
Specializing to uniform distributions on convex bodies, we show that
Theorem~\ref{thm:repi} recovers Milman's reverse Brunn-Minkowski inequality \cite{Mil86}.
One may also think of Theorem~\ref{thm:repi} as completing the usual analogy 
between the Brunn-Minkowski and entropy power inequalities (see, e.g., \cite{DCT91}).
%Theorem~\ref{thm:repi} is extended to the larger class of convex measures in \cite{BM10:repi}.

%ADDED FOR ARXIV:
This Note is in the nature of a sketch meant as an entry point for the interested reader
into the recent work by the authors in \cite{BM10:smb,BM09:maxent,BM10:repi}. Therefore, details are kept to a minimum. 
However, it is appropriate to mention that the tools relied upon in our proof of Theorem~\ref{thm:repi} and allied
results (some of which are not detailed in the sequel) include 
the localization lemma of Lov\'asz--Simonovits \cite{LS93},
reverse H\"older type inequalities such as Berwald's inequality \cite{Ber47, Bor73a}, 
the studies of the hyperplane conjecture by Bourgain \cite{Bou86}, Milman--Pajor \cite{MP89}, and Ball \cite{Bal86:phd, Bal88},
a submodularity property for the entropy of sums developed by the second-named
author \cite{Mad08:itw}, and ideas from the literature on convex measures and isoperimetry  developed
by Borell \cite{Bor74, Bor75a} and the first-named author \cite{Bob07, Bob10}.

\section{Intermediate results}
\label{sec:results}

\subsection{An equipartition property}
\label{sec:aep}

Let $X$ be a random vector taking values in $\RL^n$, and suppose its distribution has 
a density $f$ with respect to Lebesgue measure on $\RL^n$. 
The random variable 
$\h(X) = - \log f(X)$
may be thought of as the ``information content'' of $X$.
Note that the entropy is $h(X) =\E \h(X)$.

Because of the relevance of the information content in information theory, probability, and statistics, it is intrinsically 
interesting to understand its behavior. In particular, a natural question 
arises: Is it true that the information content concentrates around 
the entropy in high dimension?
In general, there is no reason for such a concentration property to hold.
% although for the very special case of Gaussians, explicit computations can be used to show that
%it does \cite{CP89}, with significant consequences for feedback capacities in information theory. 
However, the following proposition shows  that in fact, such a property
{\it holds uniformly for the entire class of log-concave densities}.  

\vspace{.15cm}
\begin{theorem}\label{prop:aep}
If $X$ has a log-concave density $f$ on $\RL^n$,
then for $0 \leq \varepsilon \leq 2$,
$$
\P\bigg\{\bigg| \frac{\h(X)}{n} - \frac{h(X)}{n}\bigg| \geq \varepsilon\bigg\} 
\leq 4 e^{-\varepsilon^2 n/16} .
$$
\end{theorem}
\vspace{.15cm}

No normalization whatsoever is required for this result, which is proved in \cite{BM10:smb}
using the localization lemma of Lov\'asz--Simonovits, and certain 
reverse H\"older type inequalities for log-concave measures. 

Equivalently, with high probability, $f(x)^{2/n}$ is very close to the entropy power
$N(X) = \exp\{\frac{2}{n}\,h(X)\}$, and
the distribution of $X$ itself is effectively the uniform distribution 
on the class of typical observables, or the ``typical set'' 
(defined to be the collection of all points $x\in\RL^n$ such that $f(x)$ 
lies between $e^{-h(X)-n\varepsilon}$ and $e^{-h(X)+n\varepsilon}$, 
for some small fixed $\varepsilon>0$).
The effective uniformity of the distribution of $X$ on some compact set, 
entailed by this concentration result, may be seen as an extension of the
{\it asymptotic equipartition property} (or  Shannon-McMillan-Breiman theorem)
to non-stationary stochastic processes with log-concave marginals (cf. \cite{BM10:smb}).

If one is more interested in the effective support rather than an effective uniformity,
one can simply consider a superlevel set (necessarily convex and compact) of the density $f$ instead of the 
annular region above. This effective support on a convex set
implied by Theorem~\ref{prop:aep} allows (see \cite{BM10:repi}) the transference of some
results from the setting of convex bodies to that of log-concave measures, in particular, 
the existence of $M$-ellipsoids  \cite{Mil86,Mil88:1,Mil88:2,Pis89:book}.
(Such a transference technique based on looking at superlevel sets
of log-concave densities has been anticipated before, e.g., by
\cite{KM05}, but  Theorem~\ref{prop:aep} refines those observations
and identifies the underlying concentration phenomenon.)

\vspace{.15cm}
\begin{corollary}\label{cor:M}
Let $\mu$ be a probability measure on $\R^n$ 
with log-concave density $f$ such that $\|f\|_{\infty} \geq 1$ (where
$\|f\|_{\infty}$ is the essential supremum and hence the maximum of $f$). 
Then there exists an ellipsoid $\mathcal{E}$ of volume 1
such that %$\mu$ may be put in a position where
$\mu(\mathcal{E})^{1/n} \geq c_M$
for some universal constant $c_M \in (0,1)$.
\end{corollary}
\vspace{.15cm}

Equivalently, for some linear volume-preserving map $u:\R^n \rightarrow \R^n$,
$\mu u^{-1}(D)^{1/n} \geq c_M$,
where $D$ is the Euclidean ball of volume one.

\subsection{Entropy and the maximal density value}
\label{sec:maxdens}

Trivially $h(X)\geq \log \|f\|_{\infty}^{-1}$. In fact, one can also bound the entropy from above using the maximal density value under
log-concavity (see \cite{BM09:maxent}).

\vspace{.15cm}
\begin{theorem}\label{prop:entcomp}
If a random vector $X$ in $\R^n$ has log-concave density $f$, then
\ben
\log\, \|f\|_{\infty}^{-1/n} \leq \nth h(X)  \leq\, 1 + \log\, \|f\|_{\infty}^{-1/n}.
\een
\end{theorem}
%\vspace{.2cm}

The hyperplane conjecture or slicing problem (cf. Bourgain \cite{Bou86} or Ball \cite{Bal88})
asserts that there exists a universal, positive constant $c$ (not depending on $n$) 
such that for any convex set $K$ of unit volume in $\mathbb{R}^n$,  there exists a hyperplane $H$
passing through its centroid such that the $(n-1)$-dimensional volume of the section
$K\cap H$ is bounded below by $c$. 
There are several equivalent formulations of the conjecture, all of a geometric or functional
analytic flavor (even the ones that nominally use probability). %Building on \cite{Bou91,Pao00},
The current best bound known, due to Klartag \cite{Kla06}, is $\Omega(n^{-1/4})$.  
Theorem~\ref{prop:entcomp} gives a purely information-theoretic formulation of the hyperplane conjecture.
For a random vector $X$ in $\RL^n$ with density $f$, 
let $D(X)$ or $D(f)$ denote its relative entropy from Gaussianity (which is the relative entropy
from the Gaussian $g$ with the same mean and covariance matrix, and also equals the difference
$h(g)-h(f)$).
The {\it Entropic Form of the Hyperplane Conjecture} asserts that for any 
log-concave density $f$ on $\RL^n$, 
$D(f)\leq cn$ for some universal constant $c$.
It is easy to see then that another equivalent form of the  hyperplane conjecture is that the entropic
distance from independence (i.e., the relative entropy of any log-concave measure from the product of its marginals)
is also bounded by $cn$ for some universal constant $c$.
As an aside, Klartag's result combined with our equivalence implies that 
$D(f)\leq \frac{1}{4} n\log n + cn$ for any log-concave $f$. This is already
the first quantitative demonstration of the spiritual closeness of log-concave measures to Gaussians,
which has been observed in qualitative ways numerous times (e.g., behavior as regards
functional inequalities). Let us note {\it en passant} that entropy plays a role in Ball's \cite{Bal03:talk} proof that
the KLS conjecture implies the hyperplane conjecture.

\vspace{-.2cm}
\section{Proof outline of Theorem~\ref{thm:repi}}
\label{sec:proofs}

The following ``submodularity'' property of the entropy functional
with respect to convolutions was obtained in \cite{Mad08:itw}:
Given independent random vectors $X$, $Y$, $Z$
in $\R^n$ with absolutely continuous distributions, we have
$$
h(X+Y+Z) + h(Z) \leq h(X+Z) + h(Y+Z)
$$
provided that all entropies are well-defined.

Let $Z\sim \text{Unif}(D)$, where $D$ is the centered Euclidean ball with volume one. 
Since $h(Z)=0$, the submodularity property implies
\ben
h(X+Y) \leq h(X+Y+Z) \leq h(X+Z)\, + \,h(Y+Z),
\een
for random vectors $X$ and $Y$ in $\R^n$ independent of each 
other and of $Z$.

Let $X$ and $Y$ have log-concave densities. Due to homogeneity of Theorem~\ref{thm:repi}, 
assume without loss of generality that $\|f\|_{\infty} \geq 1$ and 
$\|g\|_{\infty} \geq 1$. Then, our task reduces to showing that both $\calN(X+Z)$ 
and $\calN(Y+Z)$ can be bounded from above by  universal constants. 

By Corollary~\ref{cor:M}, for some affine volume preserving map 
$u:\R^n \rightarrow \R^n$, the distribution $\widetilde \mu$ of 
$\widetilde X = u(X)$ satisfies $\widetilde \mu(D)^{1/n} \geq c_M$
with a universal constant $c_M > 0$. Let $\tilde f$ denote the density of 
$\widetilde X = u(X)$. Then the density $p$ of $S = \widetilde X + Z$,
given by
$p(x) = \int_D \tilde f(x-z)\,dz = \widetilde \mu(D - x)$,
satisfies
$\|p\| \geq p(0) \geq c_M^n$.
Applying Theorem~\ref{prop:entcomp} to the random vector $S$,
$\calN(S) \leq C\, \|p\|_{\infty}^{-2/n} \leq C \cdot c_M^{-2}$,
which completes the proof.

\vspace{.2cm}
\noindent{\bf Remark 1.}
Recall C. Borell's hierarchy of convex measures on $\R^n$, classified by a parameter
$\k \in [-\infty, 1/n]$. 
%a probability measure $\mu$ on $\R^n$ is 
%called $\k$-concave, if it satisfies the Brunn-Minkowski-type inequality
%\be\label{eq:kconc-def}
%\mu \big (tA + (1-t)B \big ) \geq
%\big [\,t\mu(A)^\k + (1-t)\mu(B)^\k\big ]^{1/\k}
%\en
%for all $t \in (0,1)$ and for all Borel measurable sets $A,B \subset \R^n$ 
%with positive measure. 
In this hierarchy, $\kappa = 0$ corresponds to the class of 
log-concave measures. When $\kappa > 0$, a $\k$-concave probability measure is necessarily
compactly supported on some convex set.

%we bound the entropy $h(X)$ of a $\k$-concave random vector 
%in $\R^n$ with a positive parameter of convexity $\k$ in terms of the 
%volume of the supporting set $A$ of its distribution. 
For any random vector $X$ with values in $A$, there is a general 
upper bound $h(X) \leq \log |A|$. Using Berwald's inequality, we
provide a complementary estimate from below
depending only on the ``strength'' of convexity of
the density $f$ of $X$:
%\vspace{.3cm}
%\begin{proposition}\label{cor:cvx-ent2}
{\it Let $X$ be a random vector in $\R^n$ having an
absolutely continuous $\k$-concave distribution supported on a convex
body $A$ with $0 < \k \leq 1/n$. Then
$h(X) \geq \log |A| + n \log(\k n)$.}
%\end{proposition}
%\vspace{.3cm}
Note when $\k = 1/n$, this bound is sharp.

Assume a probability measure $\mu$ is $\k'$-concave
on $\R^n$ and a probability measure $\nu$ is $\k''$-concave on $\R^n$. 
If $\k', \k'' \in [-1,1]$ satisfy
\be\label{eq:cvx-conv}
\k' + \k'' > 0, \qquad \frac{1}{\k} = \frac{1}{\k'} + \frac{1}{\k''},
\en
then their convolution $\mu * \nu$ is $\k$-concave.
Hence, if random vectors $X_1$ and $X_2$ are independent
and uniformly distributed in convex bodies $A_1$ and $A_2$ in $\R^n$, then
the sum
$X_1 + X_2$
has a $\frac{1}{2n}$--concave distribution supported on the convex body
$A_1 +A_2$. The preceding entropy bound then implies that
$h(X_1+X_2) \geq \log |A_1+A_2| -n\log 2$.
This immediately allows one to deduce Milman's reverse 
Brunn-Minkowski inequality from Theorem~\ref{thm:repi}.

\vspace{.2cm}
\noindent{\bf Remark 2.}
Theorems~\ref{thm:repi} and \ref{prop:entcomp} have been extended
to the larger class of convex measures \cite{BM10:repi, BM09:maxent}.

%ADDED FOR ARXIV:
\vspace{.2cm}
\noindent{\bf Remark 3.}
It is instructive to write some of the above observations as entropy comparison results.
For example, the bound in Remark 1
says that if $f$ is a density of the random vector $X$ taking values in the 
convex body $A$, and if $\text{U}_A$ is the uniform distribution on $A$,
%\eqref{eq:3ub}--\eqref{eq:3lb} are equivalent to the statement
$$
0 \,\leq\, \frac{1}{n}\, D(f\|\text{U}_A) \,=\, 
\frac{1}{n}\,  h(\text{U}_A) - \frac{1}{n}\, h(X) \,\leq \, C.
$$
where $D(f\|g) =  \int f(x) \log \frac{f(x)}{g(x)}\ dx$ is the relative entropy.
In this form, the constant upper bound on the difference between relative 
entropies per coordinate expresses a dimension-free approximation of the
distribution of $X$ by the uniform distribution $\text{U}_A$. 
On the other hand, Theorem~\ref{prop:entcomp} can be rewritten
as saying that if $Z$ in $\RL^n$ is any normally distributed random vector with 
maximum density being the same as that of $X$, then 
\ben
\nth h(Z) -\half \leq\, \frac{1}{n}\, h(X)  \leq\, 
\nth h(Z) +\half.
\een
This is  a dimension-free approximation of the entropy of $X$ by a Gaussian entropy.
The entropic formulation of the hyperplane conjecture described in Section 2
is another such dimension-free approximation, both of entropy and distributional
(see \cite{BM09:maxent} for details).

\bibliographystyle{plain}
%\bibliography{$HOME/Documents-ACADEMIC/WRITINGS/CommonResources/poi,$HOME/Documents-ACADEMIC/WRITINGS/CommonResources/ik}   

\begin{thebibliography}{25}

\bibitem{Bal86:phd}
K.~Ball.
\newblock {\em Isometric problems in $\ell^p$ and sections of convex sets}.
\newblock PhD thesis, University of Cambridge, UK, 1986.

\bibitem{Bal88}
K.~Ball.
\newblock Logarithmically concave functions and sections of convex sets in
  {${\bf R}\sp n$}.
\newblock {\em Studia Math.}, 88(1):69--84, 1988.

\bibitem{Bal03:talk}
K.~M. Ball.
\newblock Information decrease along semigroups.
\newblock Talk given at conference on Banach Spaces and Convex Geometric
  Analysis, Universit{\"a}t Kiel, Germany, April 2003.

\bibitem{Ber47}
L.~Berwald.
\newblock Verallgemeinerung eines {M}ittelwertsatzes von {J}. {F}avard f\"ur
  positive konkave {F}unktionen.
\newblock {\em Acta Math.}, 79:17--37, 1947.

\bibitem{Bob07}
S.~G. Bobkov.
\newblock Large deviations and isoperimetry over convex probability measures
  with heavy tails.
\newblock {\em Electron. J. Probab.}, 12:1072--1100 (electronic), 2007.

\bibitem{Bob10}
S.~G. Bobkov.
\newblock Convex bodies and norms associated to convex measures.
\newblock {\em Probab. Theory Related Fields}, 147(1-2):303--332, 2010.

\bibitem{BM10:smb}
S.~G. Bobkov and M.~Madiman.
\newblock Concentration of the information in data with log-concave
  distributions.
\newblock {\em Ann. Probab.}, in press, {\tt arXiv:1012.5457 [math.PR]}, 2010.

\bibitem{BM09:maxent}
S.~G. Bobkov and M.~Madiman.
\newblock The entropy per coordinate of a random vector is highly constrained
  under convexity conditions.
\newblock {\tt arXiv:1006:2883 [math.PR]}, 2010.

\bibitem{BM10:repi}
S.~G. Bobkov and M.~Madiman.
\newblock Reverse {B}runn-{M}inkowski and reverse entropy power inequalities
  for convex measures.
\newblock {\em Preprint}, 2010.

\bibitem{Bor73a}
C.~Borell.
\newblock Complements of {L}yapunov's inequality.
\newblock {\em Math. Ann.}, 205:323--331, 1973.

\bibitem{Bor74}
C.~Borell.
\newblock Convex measures on locally convex spaces.
\newblock {\em Ark. Mat.}, 12:239--252, 1974.

\bibitem{Bor75a}
C.~Borell.
\newblock Convex set functions in {$d$}-space.
\newblock {\em Period. Math. Hungar.}, 6(2):111--136, 1975.

\bibitem{Bou86}
J.~Bourgain.
\newblock On high-dimensional maximal functions associated to convex bodies.
\newblock {\em Amer. J. Math.}, 108(6):1467--1476, 1986.

\bibitem{DCT91}
A.~Dembo, T.M. Cover, and J.A. Thomas.
\newblock Information-theoretic inequalities.
\newblock {\em IEEE Trans. Inform. Theory}, 37(6):1501--1518, 1991.

\bibitem{Kla06}
B.~Klartag.
\newblock On convex perturbations with a bounded isotropic constant.
\newblock {\em Geom. Funct. Anal.}, 16(6):1274--1290, 2006.

\bibitem{KM05}
B.~Klartag and V.~D. Milman.
\newblock Geometry of log-concave functions and measures.
\newblock {\em Geom. Dedicata}, 112:169--182, 2005.

\bibitem{LS93}
L.~Lov{\'a}sz and M.~Simonovits.
\newblock Random walks in a convex body and an improved volume algorithm.
\newblock {\em Random Structures Algorithms}, 4(4):359--412, 1993.

\bibitem{Mad08:itw}
M.~Madiman.
\newblock On the entropy of sums.
\newblock In {\em Proc. IEEE Inform. Theory Workshop}, pages 303--307. Porto,
  Portugal, 2008.

\bibitem{Mil86}
V.~D. Milman.
\newblock In\'egalit\'e de {B}runn-{M}inkowski inverse et applications \`a la
  th\'eorie locale des espaces norm\'es.
\newblock {\em C. R. Acad. Sci. Paris S\'er. I Math.}, 302(1):25--28, 1986.

\bibitem{Mil88:2}
V.~D. Milman.
\newblock Entropy point of view on some geometric inequalities.
\newblock {\em C. R. Acad. Sci. Paris S\'er. I Math.}, 306(14):611--615, 1988.

\bibitem{Mil88:1}
V.~D. Milman.
\newblock Isomorphic symmetrizations and geometric inequalities.
\newblock In {\em Geometric aspects of functional analysis (1986/87)}, volume
  1317 of {\em Lecture Notes in Math.}, pages 107--131. Springer, Berlin, 1988.

\bibitem{MP89}
V.~D. Milman and A.~Pajor.
\newblock Isotropic position and inertia ellipsoids and zonoids of the unit
  ball of a normed {$n$}-dimensional space.
\newblock In {\em Geometric aspects of functional analysis (1987--88)}, volume
  1376 of {\em Lecture Notes in Math.}, pages 64--104. Springer, Berlin, 1989.

\bibitem{Pis89:book}
G.~Pisier.
\newblock {\em The volume of convex bodies and {B}anach space geometry},
  volume~94 of {\em Cambridge Tracts in Mathematics}.
\newblock Cambridge University Press, Cambridge, 1989.

\bibitem{Sta59}
A.J. Stam.
\newblock Some inequalities satisfied by the quantities of information of
  {F}isher and {S}hannon.
\newblock {\em Information and Control}, 2:101--112, 1959.

\end{thebibliography}

\end{document}